# Univariate subdivision schemes for noisy data

Nira Dyn [*1], Allison Heard [†2], Kai Hormann [‡3], and Nir Sharon [§1]

[1]Department of Mathematics, Tel-Aviv University, Tel-Aviv, Israel
[2]Department of Mathematics, The University of Auckland, Auckland, New Zealand
[3]Faculty of Informatics, Università della Svizzera italiana, Lugano, Switzerland

April 17, 2013

## Abstract

We introduce and analyse univariate, linear, and stationary subdivision schemes for refining noisy data, by fitting local least squares polynomials. We first present primal schemes, based on fitting linear polynomials to the data, and study their convergence, smoothness, and basic limit functions. We provide several numerical experiments that illustrate the limit functions generated by these schemes from initial noisy data, and compare the results with approximations obtained from noisy data by an advanced local linear regression method. We conclude by discussing several extension and variants.

## 1 Introduction

In recent years, subdivision schemes have become an important tool in many applications and research areas, including animation, computer graphics, and computer aided geometric design, just to name a few [1, 16]. A subdivision scheme generates values associated with the vertices of a sequence of nested meshes, with a dense union, by repeated application of a set of local refinement rules. These rules determine the values associated with a refined mesh from the values associated with the coarser mesh. The subdivision scheme is convergent if the generated values converge uniformly to the values of a continuous function, for any set of initial values.

The particular class of interpolatory schemes consists of schemes with refinement rules that keep the values associated with the coarse mesh and only generate new values related to the additional vertices of the refined mesh. An important family of interpolatory schemes is the family of Dubuc–Deslauriers (DD) schemes [6].

Intensive studies have been carried out recently on the generalization of subdivision schemes to treat more complicated data such as manifold valued data [18, 19], matrices [17], sets [9], curves [15], and nets of functions [5]. Yet, the question how to approximate a function from its noisy samples by subdivision schemes has remained open, and it is the purpose of this paper to address this problem.

The linear and symmetric refinement rules of the DD schemes and their dual counterparts [10] are based on *local polynomial interpolation*. These schemes are stationary in the sense that the same rules are applied in each subdivision step, and their approximation order is determined by the degree of the local interpolation polynomials which are used for the refinement.

In this paper we generalize this approach and suggest linear and symmetric refinement rules based on *local polynomial approximation*, where the polynomial is determined by a least squares fit to the data. In fact, we show that the DD schemes are a special case of these schemes. The resulting schemes are stationary and well-suited for noisy input data as shown by our numerical results in Sections 3.2 and 4.4. We call these schemes least squares schemes.

[*]niradyn@post.tau.ac.il
[†]heard@math.auckland.ac.nz
[‡]kai.hormann@usi.ch
[§]nirsharo@post.tau.ac.il

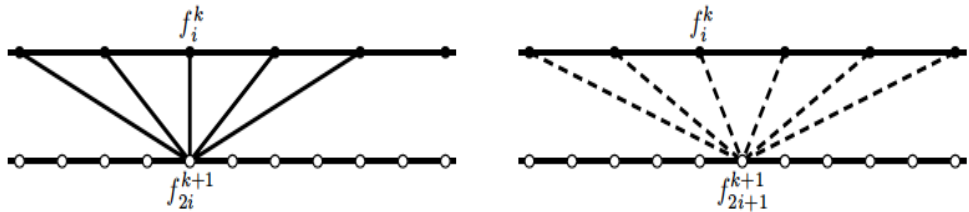

Figure 1: Schematic overview of the refinement rules in (2) for $n = 3$.

In the univariate setting that we consider, we denote by $\boldsymbol{f}^k = (f_i^k)_{i \in \mathbb{Z}}$ the data at refinement level $k \in \mathbb{N}_0$. We assume that the initial data $\boldsymbol{f}^0 = (f_i^0)_{i \in \mathbb{Z}}$ is given at the integers $\mathbb{Z}$ and that $f_i^k$ is associated with the dyadic point $t_i^k = 2^{-k} i$. The main idea of least squares subdivision is to generate the data at level $k+1$ by evaluating a polynomial that locally fits the data at level $k$ in a symmetric neighbourhood.

In particular, we use polynomials that fit the data best in the least squares sense. That is, for given data $y_1, \ldots, y_m$ at nodes $x_1, \ldots, x_m$, we are interested in the polynomial $p_d$ of some degree $d$ that minimizes the sum of squared residuals,

$$\sum_{i=1}^{m} (p_d(x_i) - y_i)^2.$$

For $m > d$ this problem has a unique solution and in Appendix A we provide a summary of the relevant theory, which also includes the case $m \leq d$. For the special case $d = 1$ and equidistant nodes $x_i = a + ih$, it turns out that the value of the linear least squares polynomial $p_1$ at the centre $c = (x_1 + \cdots + x_m)/m$ of the nodes is simply

$$p_1(c) = (y_1 + \cdots + y_m)/m. \tag{1}$$

The paper is organized as follows. We start by introducing the simplest case of least squares schemes in Section 2. These schemes are based on primal refinement rules and on best fitting linear polynomials to symmetric data. This is a one parameter family of schemes, with the number of data points as the parameter. We prove convergence and smoothness of these schemes and investigate properties of the corresponding basic limit functions. The construction of least squares schemes based on best fitting polynomials of higher degrees and on dual refinement rules is postponed to Section 4. In Section 3 we review a statistical model for fitting noisy data, analyse the suitability of the primal, least squares schemes of degree 1 for dealing with this kind of data, and provide several numerical examples. Further results for primal schemes based on best fitting polynomials of higher order are presented in Section 4.4. Throughout this paper we use several well-known properties of least squares polynomials. A short survey of these properties of least squares polynomials and a method for the efficient evaluation of our schemes are given in Appendix A.

## 2  Primal least squares schemes of degree 1

We start by considering the simplest least squares subdivision scheme $S_n$ for $n \geq 1$, which generates the data at level $k + 1$ as follows. On the one hand, the value $f_{2i}^{k+1}$, which replaces $f_i^k$, is determined by fitting a linear polynomial to the $2n - 1$ data values in a symmetric neighbourhood around $t_i^k$ and evaluating it at the associated dyadic point $t_i^k = t_{2i}^{k+1}$. On the other, the scheme computes the new value $f_{2i+1}^{k+1}$ between $f_i^k$ and $f_{i+1}^k$ by evaluating the linear least squares polynomial with respect to the nearest $2n$ data values halfway between the corresponding dyadic points $t_i^k$ and $t_{i+1}^k$, namely at $t_{2i+1}^{k+1}$. In this construction the parameter $n$ controls the locality of the scheme and we study its effect in Section 3.

Following (1), the refinement rules of $S_n$ turn out to be

$$f_{2i}^{k+1} = \frac{1}{2n-1} \sum_{j=-n+1}^{n-1} f_{i+j}^k \quad \text{and} \quad f_{2i+1}^{k+1} = \frac{1}{2n} \sum_{j=-n+1}^{n} f_{i+j}^k. \tag{2}$$

Consequently, the symbol [7] of the scheme is

$$a_n(z) = \frac{1}{2n} \sum_{j=-n+1}^{n} z^{2j-1} + \frac{1}{2n-1} \sum_{j=-n+1}^{n-1} z^{2j}. \tag{3}$$

It follows from the symmetry of the points determining the linear least squares polynomials, that $a_n(z) = a_n(1/z)$, hence the scheme is odd symmetric [11]. As the data at level $k+1$ depends on at most $2n$ values at level $k$, we conclude that $S_n$ is a primal $2n$-point scheme. An example of the refinement dependencies for $S_3$ is shown in Figure 1, and the masks of the first three schemes are

$$\begin{aligned} \boldsymbol{a}_1 &= & [1,2,1] & \quad /\, 2, \\ \boldsymbol{a}_2 &= & [3,4,3,4,3,4,3] & \quad /\, 12, \\ \boldsymbol{a}_3 &= [5,6,5,6,5,6,5,6,5,6,5] & & /\, 30, \end{aligned}$$

Note that the scheme $S_1$ is the interpolating 2-point scheme, which generates piecewise linear functions in the limit.

## 2.1 Convergence and smoothness

Following the usual definition of convergence in [8, Chapter 2], we denote the limit of a convergent subdivision scheme $S$ for initial data $\boldsymbol{f}^0$ by $S^\infty \boldsymbol{f}^0$.

The explicit form of the symbol in (3) implies that $a_n(1) = 2$ and $a_n(-1) = 0$, which are necessary conditions for $S_n$ to be convergent [7, Proposition 2.1]. Following the analysis in [7], we define

$$q_n(z) = \frac{a_n(z)}{1+z} = \frac{1}{2n(2n-1)} \left( \sum_{j=-n+1}^{n-1} (n-j) z^{2j-1} + \sum_{j=-n+1}^{n-1} (n+j) z^{2j} \right) \tag{4}$$

and conclude the convergence of $S_n$ by analysing the norm of the subdivision scheme with symbol $q_n$.

**Theorem 1.** *The least squares subdivision scheme $S_n$ is convergent for $n \geq 1$.*

*Proof.* It follows from (4) that the norm of the subdivision scheme with symbol $q_n$ is

$$\|S_{[q_n]}\|_\infty = \max\left\{ \frac{1}{2n(2n-1)} \sum_{j=-n+1}^{n-1} |n-j|, \frac{1}{2n(2n-1)} \sum_{j=-n+1}^{n-1} |n+j| \right\}$$
$$= \frac{1}{2n(2n-1)} \sum_{j=1}^{2n-1} j = \frac{1}{2}.$$

According to [7, Theorem 3.2], the scheme $S_n$ is therefore convergent. $\square$

Note that the norm of the scheme $S_{[q_n]}$ is the least possible value as is the case for the uniform B-spline schemes, indicating "quickest" possible convergence. The structure of $q_n$ further reveals that the limit functions generated by $S_n$ are $C^1$.

**Theorem 2.** *The least squares subdivision scheme $S_n$ generates $C^1$ limit functions for $n \geq 2$.*

| $n$ | 2 | 3 | 4 | 5 | 6 | 7 | 8 | 9 | 10 |
|---|---|---|---|---|---|---|---|---|---|
| $\rho_n$ | 1.649 | 1.777 | 1.816 | 1.794 | 1.786 | 1.776 | 1.771 | 1.761 | 1.753 |

Table 1: Lower bounds $\rho_n$ on the Hölder regularity of the schemes $S_n$.

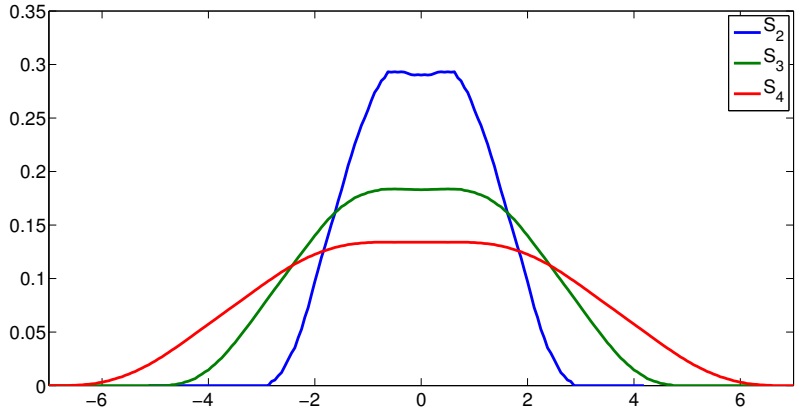

Figure 2: Basic limit functions of the schemes $S_2$, $S_3$, and $S_4$.

*Proof.* It is known [7, Theorems 3.2 and 3.4] that in order to prove the theorem, it is sufficient to show that the scheme with symbol $2q_n$ is convergent. By (4),

$$2q_n(1) = 2 \qquad \text{and} \qquad 2q_n(-1) = 0,$$

hence $S_{[2q_n]}$ satisfies the necessary conditions for convergence. Since the coefficients of the symbol are all positive and there are at least three such coefficients for $n \geq 2$, it follows from [3, Theorem 3.3] that the scheme is convergent, and so $S_n$ generates $C^1$ limit functions. □

The statement in Theorem 2 is confirmed by the numerical results presented in Table 1, which were obtained by using 16 iterations of the algorithm in [12] to compute lower bounds on the Hölder regularity.

## 2.2 The basic limit function

Let us denote by $\boldsymbol{\delta}$ the sequence which is zero everywhere except at 0, where it is 1. The *basic limit function* of the convergent subdivision scheme $S_n$ is then defined as

$$\phi_n = S_n^\infty \boldsymbol{\delta}. \tag{5}$$

Some examples of $\phi_n$ for small values of $n$ are shown in Figure 2.

Many properties of a linear subdivision scheme can be derived from its basic limit function. In particular, due to linearity, the limit function generated from the initial data $\boldsymbol{f}^0 = (f_i^0)_{i \in \mathbb{Z}}$ by the scheme $S_n$ has the form

$$(S_n^\infty \boldsymbol{f}^0)(x) = \sum_{j \in \mathbb{Z}} f_j^0 \phi_n(x - j). \tag{6}$$

Our first observation is that the support of $\phi_n$ is $[-2n+1, 2n-1]$, because $S_n$ is a primal $2n$-point scheme [6]. Moreover, $\phi_n$ is positive inside its support, because the coefficients of the mask $\boldsymbol{a}_n$ are positive in the mask's support, and $\phi_n$ has the partition of unity property

$$\sum_{j \in \mathbb{Z}} \phi_n(x - j) = 1, \tag{7}$$

due to the reproduction of constant polynomials by $S_n$.

The simple structure of $\boldsymbol{a}_n$ further allows us to derive several interesting properties regarding the values of the basic limit function $\phi_n$ at the integers. These values are of importance, because they constitute the *filter* which operates on the initial data and generates the final values at the integers. Taking into account that $\phi_n$ is continuous and therefore vanishes at the end points of its support, we conclude from (6) that the limit at the integers $k \in \mathbb{Z}$ is

$$(S_n^\infty \boldsymbol{f}^0)(k) = \sum_{j=-2n+2}^{2n-2} f_{k-j}^0 \phi_n(j). \tag{8}$$

The non-zero values of $\phi_n$ at the integers constitute an eigenvector $\boldsymbol{v} = \big(\phi_n(-2n+2), \ldots, \phi_n(2n-2)\big)$ corresponding to the eigenvalue 1 of the transposed subdivision matrix [7], which in this case is the $(4n-3) \times (4n-3)$ column stochastic, two-slanted band matrix

$$A_n = \begin{pmatrix} r & s & 0 & 0 & 0 & & 0 & 0 & 0 \\ r & s & r & s & 0 & \cdots & 0 & 0 & 0 \\ r & s & r & s & r & & 0 & 0 & 0 \\ & & \vdots & & & \ddots & & \vdots & \\ r & s & r & s & r & & r & s & 0 \\ r & s & r & s & r & \cdots & r & s & r \\ 0 & s & r & s & r & & r & s & r \\ & & \vdots & & & \ddots & & \vdots & \\ 0 & 0 & 0 & 0 & 0 & \cdots & 0 & s & r \end{pmatrix}$$

with entries $r = 1/(2n-1)$ and $s = 1/(2n)$.

The odd symmetry of the mask $\boldsymbol{a}_n$ guarantees that $\phi_n$ is a symmetric function. Thus, the eigenvector $\boldsymbol{v}$ is also symmetric, as indicated by the structure of $A_n$. Taking these symmetries into account, we get that the vector $\tilde{\boldsymbol{v}} = \big(\phi_n(-2n+2), \ldots, \phi_n(0)\big)$ is an eigenvector corresponding to the eigenvalue 1 of the $(2n-1) \times (2n-1)$ matrix

$$\tilde{A}_n = \begin{pmatrix} r & s & 0 & 0 & 0 & & 0 & 0 & 0 & 0 & 0 \\ r & s & r & s & 0 & \cdots & 0 & 0 & 0 & 0 & 0 \\ r & s & r & s & r & & 0 & 0 & 0 & 0 & 0 \\ & & \vdots & & & \ddots & & & \vdots & & \\ r & s & r & s & r & & r & s & 0 & 0 & 0 \\ r & s & r & s & r & & r & s & r & s & 0 \\ r & s & r & s & r & \cdots & r & s & r & 2s & r \\ r & s & r & s & r & & r & 2s & 2r & 2s & r \\ r & s & r & s & r & & 2r & 2s & 2r & 2s & r \\ & & \vdots & & & \ddots & & & \vdots & & \\ r & s & r & 2s & 2r & & 2r & 2s & 2r & 2s & r \\ r & 2s & 2r & 2s & 2r & \cdots & 2r & 2s & 2r & 2s & r \\ 2r & 2s & 2r & 2s & 2r & & 2r & 2s & 2r & 2s & r \end{pmatrix}.$$

The particular structure of $\tilde{A}_n$ allows us to derive the following observation.

**Lemma 3.** *The values of $\phi_n$ at the non-positive integers in its support are strictly increasing,*

$$0 < \phi_n(-2n+2) < \phi_n(-2n+3) < \cdots < \phi_n(-1) < \phi_n(0).$$

*Moreover,*

$$\phi_n(-n) = \frac{n-1}{2n-1} \phi_n(0).$$

*Proof.* Note that each row of $\tilde{A}_n$ is equal to the previous row plus at least one positive terms. Since $\tilde{\boldsymbol{v}}$ satisfies $\tilde{A}_n \tilde{\boldsymbol{v}} = \tilde{\boldsymbol{v}}$ and its components $\tilde{v}_i = \phi_n(i - 2n + 1)$, $i = 1, \ldots, 2n-1$, are positive, the latter must be strictly increasing.

To establish the second statement consider the $(n-1)$-th and the last rows of $\tilde{A}_n$,

$$\tilde{\boldsymbol{\alpha}}_{n-1} = (r, s, r, s, \ldots, r, s, 0) \quad \text{and} \quad \tilde{\boldsymbol{\alpha}}_{2n-1} = (2r, 2s, 2r, 2s, \ldots, 2r, 2s, r),$$

and note that

$$\tilde{\boldsymbol{\alpha}}_{2n-1} = 2\tilde{\boldsymbol{\alpha}}_{n-1} + (0, 0, \ldots, 0, r).$$

Then, since

$$\tilde{v}_{n-1} = \tilde{\boldsymbol{\alpha}}_{n-1} \tilde{\boldsymbol{v}}$$

and

$$\tilde{v}_{2n-1} = \tilde{\boldsymbol{\alpha}}_{2n-1} \tilde{\boldsymbol{v}} = 2\tilde{\boldsymbol{\alpha}}_{n-1} \tilde{\boldsymbol{v}} + r\tilde{v}_{2n-1} = 2\tilde{v}_{n-1} + r\tilde{v}_{2n-1},$$

the second statement follows directly from the definition of $\tilde{\boldsymbol{v}}$ as $r = \frac{1}{2n-1}$. □

By the symmetry of $\phi_n$, the statements of Theorem 3 hold analogously for the values of $\phi_n$ at the non-negative integers. We are now ready to turn this result, in view of (8) into a property regarding the values of the limit functions of our scheme at the integers.

**Corollary 4.** *The least squares subdivision scheme $S_n$ acts as a filter on the initial data $\boldsymbol{f}^0$, such that $(S^\infty \boldsymbol{f}^0)(k)$ is a convex combination of $f^0_{k-j}$, $|j| \leq 2n-2$. The weights corresponding to $f^0_{k-j}$ for $|j| \geq n$ are at most 1/2 the maximal weight corresponding to $f^0_k$.*

Lemma 3 clarifies the behaviour of $\phi_n$ over the integers. However, as seen in Figure 2, this appear to be also true for the values between the integers, and it motivates to further analyse $\phi_n$. Next, we extract several bounds on the basic limit function.

The masks corresponding to the refinement rules (2) are positive. Thus, for a non-negative data such as $\delta$ we have $\left\|S_n^{k_1}(\delta)\right\|_\infty \leq \left\|S_n^{k_2}(\delta)\right\|_\infty$ for any integers $k_1 > k_2 > 0$. In other words, we can bound $\|\phi_n\|_\infty$, for example, by

$$\|\phi_n\|_\infty = \|S_n^\infty(\delta)\|_\infty \leq \left\|S_n^1(\delta)\right\|_\infty = \frac{1}{(2n-1)}. \tag{9}$$

The bound (9) suggests that the basic limit functions are asymptotically tends to zero as $n$ grows.

Our numerical computations indicate that the behaviour of $\phi_n$ between consecutive integers is not far from being monotone, and as $n$ grows even becomes constant. To see this, we study the asymptotic behaviour of the derivative $\phi'_n$ which exists since $S_n$ generates $C^1$ limits (Theorem 2). By (4) and (5) we have the relation [7, Section 2.3]

$$\phi'_n(x) = S^\infty_{[2q_n]}(\Delta\delta)(x),$$

where $\Delta$ is the forward differences operator, namely $(\Delta\delta)_0 = -1$, $(\Delta\delta)_{-1} = 1$, and zero otherwise. By (4)

$$(2q_n)_{2j-1} = \frac{2}{2n(2n-1)}(n-j), \quad j = -n+1, \ldots, n-1$$

$$(2q_n)_{2j} = \frac{2}{2n(2n-1)}(n+j), \quad j = -n+1, \ldots, n-1.$$

The scheme $S_{[2q_n]}$ has a positive mask, which entails $\|\phi'_n\|_\infty \leq \|S_{[2q_n]}\Delta\delta\|_\infty$. A direct calculation yields $\left(S_{[2q_n]}\Delta\delta\right)_j = (2q_n)_{2j+1} - (2q_n)_{2j}$. Therefore,

$$\|S_{[2q_n]}\Delta\delta\|_\infty = \max_j\{|(2q_n)_{2j+1} - (2q_n)_{2j}|\} = 2\max\{\frac{2n+3}{2n(2n-1)}, \frac{1}{2n}\} = \frac{2n+3}{2n-1}\frac{1}{n}. \tag{10}$$

To conclude, we obtain an asymptotic bound for the derivative of $\phi_n$

$$\|\phi_n'\|_\infty \sim 1/n,$$

suggests that $\phi_n$ becomes almost constant between the integers. This observation supports our simulations and when combined with Lemma 3 and (9) enhance our understanding of the basic limit function $\phi_n$.

## 3 The schemes applied to noisy data

The schemes $\{S_n\}_{n>1}$ are designed to deal with noisy data, which is supported by the following discussions and experiments. We first introduce a statistical model and then compare the performance of our schemes and an advanced local linear regression method.

### 3.1 Statistical considerations

Let $f\colon \mathbb{R} \to \mathbb{R}$ be a continuous scalar function and suppose we are given a discrete set of noisy samples

$$y_i = f_i + \varepsilon_i, \quad i \in \mathbb{Z}, \tag{11}$$

where $f_i = f(ih)$ and $\varepsilon_i$ is a normally distributed error with zero mean and variance $\sigma^2$. As an estimator of $f$ we use the limit (6) of $S_n$, that is

$$\hat{f}(x) = \sum_{j \in \mathbb{Z}} y_j \phi_n(x - j). \tag{12}$$

Note that $\hat{f}(x)$ is a random variable and the estimation quality of $\hat{f}$ is given by the expectation of the squared error.

**Lemma 5.** *With* $\mathrm{E}$ *denoting the expectation operator, the expected squared error for $x \in \mathbb{R}$ is*

$$\mathrm{E}\big[(\hat{f}(x) - f(x))^2\big] = \sigma^2 \sum_{j \in \mathbb{Z}} \phi_n(x - j)^2 + \left(\sum_{j \in \mathbb{Z}} f(j)\phi_n(x - j) - f(x)\right)^2 \tag{13}$$

*Proof.* By (12) we have

$$\operatorname{Var} \hat{f}(x) = \sum_{i \in \mathbb{Z}} \operatorname{Var}(y_i) \phi_n(x - i)^2 = \sigma^2 \sum_{i \in \mathbb{Z}} \phi_n(x - i)^2$$

and

$$\mathrm{E}[\hat{f}(x)] = \sum_{i \in \mathbb{Z}} \mathrm{E}(y_i) \phi_n(x - i) = \sum_{i \in \mathbb{Z}} f_i \, \phi_n(x - i).$$

Since $\operatorname{Var} \hat{f}(x) = \mathrm{E}\big[\hat{f}(x)^2\big] - \mathrm{E}[\hat{f}(x)]^2$, Equation (13) follows from

$$\begin{aligned}
\mathrm{E}\big[(\hat{f}(x) - f(x))^2\big] &= \mathrm{E}\big[\hat{f}(x)^2 - 2\hat{f}(x)f(x) + f(x)^2\big] \\
&= \mathrm{E}\big[\hat{f}(x)^2\big] - 2\,\mathrm{E}[\hat{f}(x)]f(x) + f(x)^2 \\
&= \operatorname{Var} \hat{f}(x) + \mathrm{E}[\hat{f}(x)]^2 - 2\,\mathrm{E}[\hat{f}(x)]f(x) + f(x)^2 \\
&= \operatorname{Var} \hat{f}(x) + \big(\mathrm{E}[\hat{f}(x)] - f(x)\big)^2.
\end{aligned}$$

□

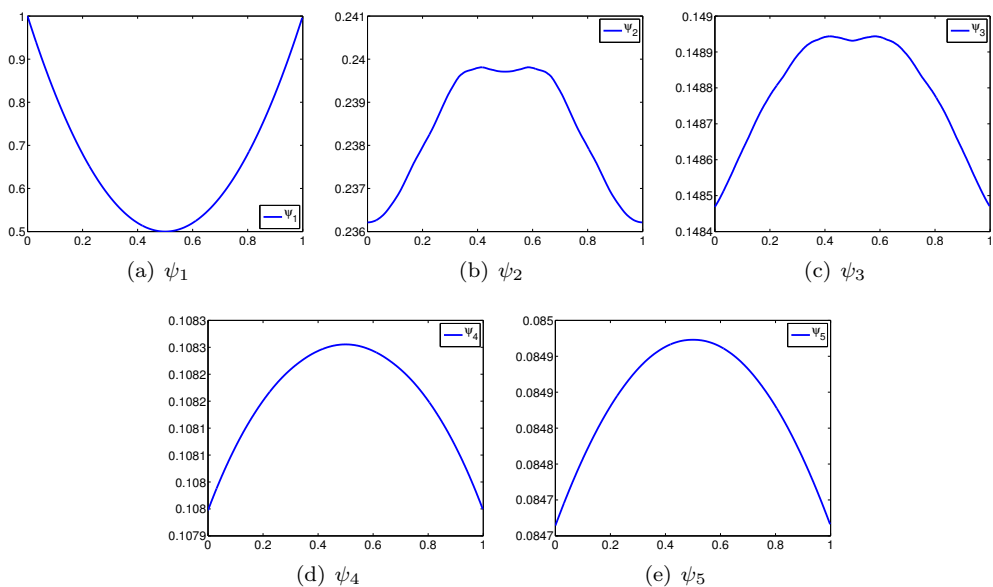

Figure 3: $\psi_n$ for the values $n = 1, \ldots, 5$. Note the different scale in each sub-figure.

By Lemma 5 the expectation of the squared error consists of two terms. The first term is the product of the variance of the noise $\sigma^2$ and the function

$$\psi_n(x) = \sum_{i \in \mathbb{Z}} \phi_n(x-i)^2. \tag{14}$$

The second term is the square of the deterministic approximation error corresponding to data without noise. We first study $\psi_n$ and come back to the second term later.

Figure 3 presents several numerical evaluations of $\psi_n$. Note how $\psi_n$ becomes almost constant as $n$ grows. We further analyse $\psi_n$ in order to establish this observation. By (7) and the positivity of $\phi_n$ we have

$$\psi_n = \sum_{j \in \mathbb{Z}} \phi_n(x-j)^2 \leq 1, \tag{15}$$

with a strict inequality for $n > 1$, namely for non-interpolatory schemes. Furthermore, as can be seen in Figure 3(a), for the interpolation scheme we have $\psi_1(x) = 1$, $x \in \mathbb{Z}$. It is a common knowledge that interpolation is not appropriate for noisy data. Indeed, (15) explains that the effect of the noise on the expected square error is bigger in the interpolatory scheme based on interpolation by local linear polynomials than in schemes based on local least squares linear polynomials. In other words, a small $\|\psi_n\|_\infty$ guarantees less effect of the noise on our estimator.

In the following we summarize the properties of $\psi_n$ and derive a bound on $\|\psi_n\|_\infty$.

**Theorem 6.** *The function $\psi_n$ in (14) is positive, symmetric, and periodic with period 1. Moreover,*

$$\int_0^1 \psi_n(x)dx \sim \frac{1}{n},$$

*and*

$$\|\psi_n\|_\infty \leq C\frac{1}{n},$$

*with a real constant $C$.*

*Proof.* By definition, $\psi_n$ is positive and periodic, and since $\phi_n$ is compactly supported, $\psi_n$ is finite. The symmetry of $\phi_n$ implies the symmetry of $\psi_n$, and it follows from the periodicity that $\psi_n$ is also symmetric about $1/2$. In addition,

$$\|\phi_n\|_2^2 = \int_{\mathbb{R}} \phi_n(x)^2 dx = \sum_{j \in \mathbb{Z}} \int_0^1 \phi_n(x-j)^2 dx = \int_0^1 \sum_{j \in \mathbb{Z}} \phi_n(x-j)^2 dx = \int_0^1 \psi_n(x) dx,$$

where the interchange of summation and integration is valid because the sum is finite. Using the bound (9) we have

$$\|\phi_n\|_2^2 \leq \frac{4n-2}{(2n-1)^2} = \frac{2}{2n-1}.$$

The latter combined with (10) and the general observation that a $C^1$ function $f$ on a closed segment $I$ satisfies $\max_I(f) \leq |I| \max_I(f') + \min_I(f)$, yield

$$\|\psi_n\|_\infty \leq \frac{2n+1}{n(2n-1)} + \frac{2}{2n-1} = \frac{4n+1}{2n^2-n} \leq 2\frac{1}{n} + \mathcal{O}(\frac{1}{n^2}).$$

We used the trivial fact that $\min_{[0,1]}(\psi_n) \leq \int_0^1 \psi_n(x) dx$. □

We note that for the proof of Theorem 6 we use the bound (9), established with a single refinement of the data $\boldsymbol{\delta}$. Better bounds can be obtained using a few more refinements, however these bounds are merely a slight improvement for the constant of the leading term $1/n$ of $\|\psi_n\|_\infty$.

The second term of the expected square error in (13) is the deterministic error or the approximation error. We use the approximation order as a standard measure for the quality of the approximation, see e.g., [12, Chapter 7]. For the case of schemes based on linear least squares polynomials, the approximation order is $h^2$, where $h$ is the distance between the sampled points of the initial data. This observation follows from the polynomial reproduction property of our schemes, that is the reconstruction of any linear polynomial from its samples.

In conclusion, there is a trade-off between the deterministic approximation error and the effect of the noise on the expected square error. In particular, higher values of $n$ decrease the effect of the noise but increase the deterministic error due to averaging of the values $\{f(i)\}$ by weights with a large support.

## 3.2 Numerical simulations

We illustrate the performance of some of the schemes by several numerical experiments, starting from noisy data. We compare their performance with today's state-of-the-art algorithm for local fitting of noisy data, namely with local linear regression (LLR).

This local estimator around a given data point $x^*$ is obtained by including kernel weights into the least squares minimization problem in the neighbourhood of $x^*$,

$$\min_{\alpha, \beta} \sum_{i=0}^n (y_i - \alpha - \beta(x_i - x^*))^2 \operatorname{Ker}(x_i - x^*).$$

This approach can be generalized to higher degree polynomials as well. For more details see [14, Chapter 4]. Although the concept of LLR is rather simple, it is one of the most important statistical approaches used.

We take the LLR, which is based on the normal kernel with the parameters of the kernel chosen dynamically, and compare it with the limits of the subdivision schemes $S_3$ and $S_5$. In the first examples, which are presented in Figures 4 and 5, the sampled function is the slowly varying function $f(x) = \sin(x/10) + (\frac{x}{50})^2$. This function is plotted in both figures, as well as the noisy data $\{y_i\}$ created with a normally distributed noise. As discussed in Subsection 3.1, the subdivision scheme $S_3$ with the smaller support, is more sensitive to the variance of the noise than $S_5$. The latter generates a limit which almost coincides with the results of the LLR, while the limit of $S_3$ is inferior.

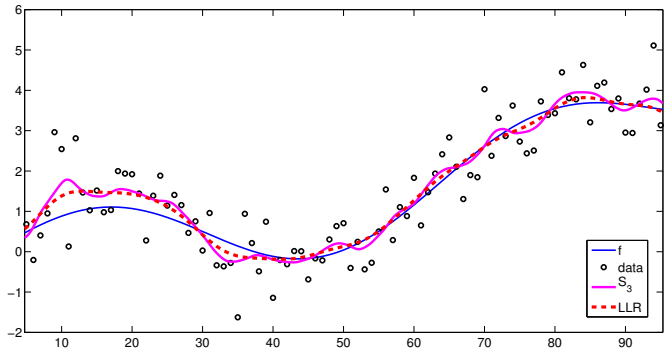

Figure 4: Comparison of the performance of $S_3$ and the local linear regression (LLR).

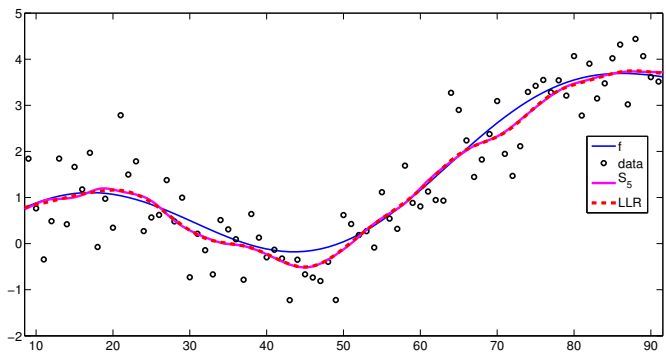

Figure 5: Comparison of the performance of $S_5$ and the local linear regression (LLR).

In Figure 6 we present a different comparison between $S_3$ and the LLR for data taken from the oscillatory function $f(x) = \cos(0.4x) + (\frac{x}{40} - 1)^3$, with high level of noise. For this data the limit of the subdivision scheme is superior.

In Figure 7 we provide a comparison between $S_5$ and the LLR where the data is sampled from the step function $f(x) = \begin{cases} 1, & x \geq 50 \\ 0, & \text{otherwise} \end{cases}$, with high level of noise. Similar to the previous example, for this data the limit of the subdivision scheme seems to give a better fit.

We observe that the numerical simulations support our understanding about the trade-off between the effect of the noise and the deterministic approximation error, as discussed in Section 3.1.

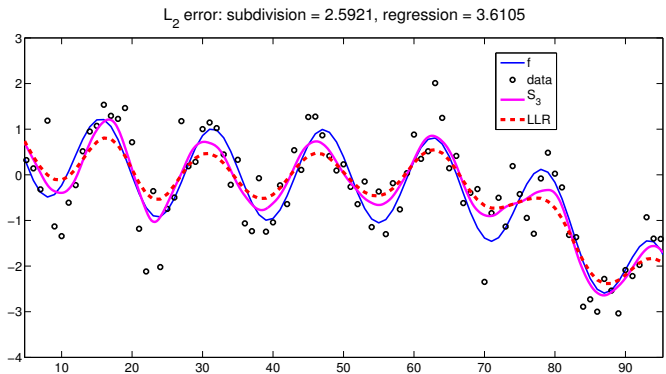

Figure 6: Comparison of the performance of $S_3$ and the local linear regression (LLR) in case of an oscillatory function with high noise level. The $L_2$ errors of the estimators are displayed above.

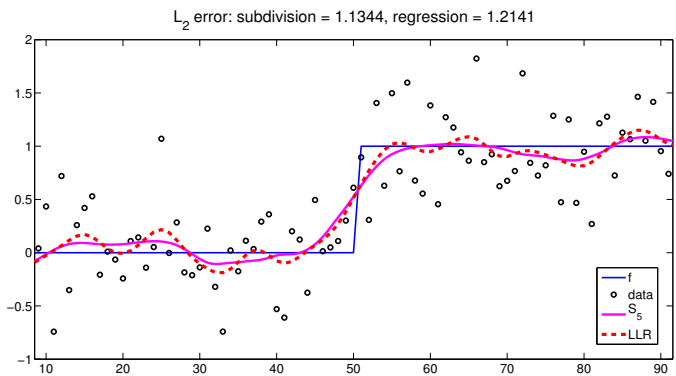

Figure 7: Comparison of the performance of $S_3$ and the local linear regression (LLR) in case of an oscillatory function with high noise level. The $L_2$ errors of the estimators are displayed above.

## 4 Extensions and variants

The family of primal least squares schemes of degree 1 can be extended in several ways. We first discuss the extension to dual schemes (Section 4.1), as well as minor variations of both primal and dual schemes (Section 4.2). A further extension relies on fitting least squares polynomials of higher degree (Section 4.3) and we provide a few numerical examples of such schemes (Section 4.4).

### 4.1 Dual least squares schemes of degree 1

While the idea of the schemes $S_n$ in Section 2 was to fit linear least squares polynomials and to evaluate them in a *primal* way, that is, at the points and the midpoints of the mesh at level $k$, another option is to design subdivision schemes based on *dual* evaluation [11]. The dual least squares scheme $\bar{S}_n$ is obtained by fitting a linear polynomial to the $2n$ data values at the points $t_{i-n+1}^k, \ldots, t_{i+n}^k$ at level $k$ and then evaluating this polynomial at $1/4$ and $3/4$ between $t_i^k$ and $t_{i+1}^k$ to compute the new data $f_{2i}^{k+1}$ and $f_{2i+1}^{k+1}$. Figure 8 shows an example of the refinement dependencies for $\bar{S}_3$.

The refinement rules of the dual schemes are slightly more complicated to derive than those of the

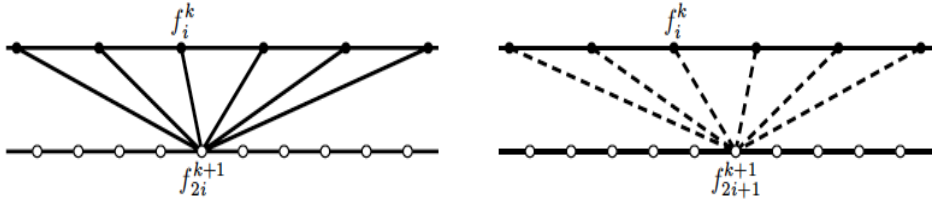

Figure 8: Schematic overview of the dual refinement rules in (16) for $n = 3$.

| $n$ | 2 | 3 | 4 | 5 | 6 | 7 | 8 | 9 | 10 |
|---|---|---|---|---|---|---|---|---|---|
| $\bar{\rho}_n$ | 2.285 | 2.647 | 2.729 | 2.677 | 2.664 | 2.633 | 2.616 | 2.594 | 2.577 |

Table 2: Lower bounds $\bar{\rho}_n$ on the Hölder regularity of the schemes $\bar{S}_n$.

primal schemes, but they still have a rather simple closed form. The refinement rules of $\bar{S}_n$ are

$$f_{2i}^{k+1} = \frac{1}{2n} \sum_{j=-n+1}^{n} \left(1 - \frac{6j-3}{8n^2-2}\right) f_{i+j}^k \quad \text{and} \quad f_{2i+1}^{k+1} = \frac{1}{2n} \sum_{j=-n+1}^{n} \left(1 + \frac{6j-3}{8n^2-2}\right) f_{i+j}^k. \quad (16)$$

The corresponding symbol is

$$\bar{a}_n(z) = \frac{1}{2n} \sum_{j=-n}^{n-1} \left(1 + z + \frac{6j+3}{8n^2-2}(1-z)\right) z^{2j}, \quad (17)$$

and it is easy to verify that $\bar{a}_n(z)z = \bar{a}_n(1/z)$, which confirms $\bar{S}_n$ to be an even symmetric scheme [11]. Overall we conclude that $\bar{S}_n$ is a dual $2n$-point scheme and the support of its basic limit function $\bar{\phi}_n$ is $[-2n, 2n-1]$. The masks of the first three schemes are

$$\bar{a}_1 = \qquad\qquad [1, 3, 3, 1] \qquad\qquad / 4,$$
$$\bar{a}_2 = \qquad [7, 13, 9, 11, 11, 9, 13, 7] \qquad / 40,$$
$$\bar{a}_3 = [55, 85, 61, 79, 67, 73, 73, 67, 79, 61, 85, 55] / 420.$$

We recognize $\bar{S}_1$ as Chaikin's corner cutting scheme [4].

The proofs of Theorems 1 and 2 carry over to the dual schemes, and so the limit functions generated by $\bar{S}_n$ are at least $C^1$ for $n \geq 1$. However, unlike the primal schemes, the symbols of the dual schemes are divisible by $(1+z)^3$, and so they may potentially generate $C^2$ limits. But there is no simple proof as for $C^1$ in Theorem 2, because the symbol $4\bar{a}_n(z)/(1+z)^2$ has negative coefficients. Table 2 lists lower bounds on the Hölder regularity of the first few schemes, computed using 16 iterations of the algorithm in [12] and demonstrates that the limits of $\bar{S}_n$ are in fact $C^2$, at least for $2 \leq n \leq 10$.

### 4.2 Variants of linear least squares schemes

In addition to the dual $2n$-point schemes $\bar{S}_n$, it is also possible to define dual $(2n+1)$-point schemes. These schemes fit a linear polynomial to the $2n+1$ data values in a symmetric neighbourhood around $f_i^k$ and evaluate it at $1/4$ the distance to the left and right neighbours to define the new data $f_{2i-1}^{k+1}$ and $f_{2i}^{k+1}$. The resulting refinement rules are

$$f_{2i-1}^{k+1} = \frac{1}{2n+1} \sum_{j=-n}^{n} \left(1 - \frac{3j}{4n(n+1)}\right) f_{i+j}^k \quad \text{and} \quad f_{2i}^{k+1} = \frac{1}{2n+1} \sum_{j=-n}^{n} \left(1 + \frac{3j}{4n(n+1)}\right) f_{i+j}^k,$$

and the support of the corresponding basic limit function is $[-2n-1, 2n]$. The masks of the first three schemes of this kind are

$$n = 1: \qquad [5, 11, 8, 8, 11, 5] \qquad /\, 24,$$
$$n = 2: \qquad [6, 10, 7, 9, 8, 8, 9, 7, 10, 6] \qquad /\, 40,$$
$$n = 3: [13, 19, 14, 18, 15, 17, 16, 16, 17, 15, 18, 14, 19, 13]\, /\, 112.$$

Similarly, we can define primal $(2n+1)$-point schemes as variants of the primal $2n$-point schemes $S_n$. We simply replace the refinement rule for $f_{2i}^{k+1}$ in (2) by

$$f_{2i}^{k+1} = \frac{1}{2n+1} \sum_{j=-n}^{n} f_{i+j}^{k}$$

and keep the rule for $f_{2i+1}^{k+1}$. For these schemes, the support of the basic limit function is $[-2n, 2n]$, and the masks of the first three schemes are

$$n = 1: \qquad [2, 3, 2, 3, 2] \qquad /\, 6,$$
$$n = 2: \qquad [4, 5, 4, 5, 4, 5, 4, 5, 4] \qquad /\, 20,$$
$$n = 3: [6, 7, 6, 7, 6, 7, 6, 7, 6, 7, 6, 7, 6]\, /\, 42.$$

Adapting the proofs of Theorems 1 and 2, one can show that both variants generate $C^1$ limit functions, and our numerical results demonstrate that the dual $(2n+1)$-point schemes are even $C^2$ for $1 \le n \le 10$.

## 4.3 Least squares schemes of higher degree

The least squares schemes of degree 1 reproduce linear polynomials by construction, but they do not reproduce polynomials of higher degree, and so their approximation order is only $O(h^2)$, unless the data is being pre-processed [11]. We can improve this by using least squares polynomials of some higher degree $d > 1$. To this end, let $p_{n,i}^d$ be the least squares polynomial of degree $d$ for the $2n-1$ data

$$(t_{i+j}^k, f_{i+j}^k), \qquad j = -n+1, \ldots, n-1$$

in a symmetric neighbourhood of $t_{2i}^{k+1}$, and let $\tilde{p}_{n,i}^d$ be the polynomial of degree $d$ that fits the $2n$ data

$$(t_{i+j}^k, f_{i+j}^k), \qquad j = -n+1, \ldots, n$$

in a symmetric neighbourhood of $t_{2i+1}^{k+1}$. The polynomials $p_{n,i}^k$ and $\tilde{p}_{n,i}^k$ are well-defined for $d < 2n-1$ and $d < 2n$, respectively; see Appendix A.1 for details.

The primal $2n$-point least squares scheme of degree $d$ is then characterized by the refinement rules

$$f_{2i}^{k+1} = p_{n,i}^d(t_i^k) \qquad \text{and} \qquad f_{2i+1}^{k+1} = \tilde{p}_{n,i}^d\bigl((t_i^k + t_{i+1}^k)/2\bigr), \qquad (18)$$

which simplifies to the rules in (2) for $d = 1$. The resulting subdivision scheme $S_n^d$ reproduces polynomials of degree $d$ by construction and thus has approximation order $O(h^{d+1})$. It is well-defined for $d < 2n$, even though for $d = 2n-1$ the rule for $f_{2i}^{k+1}$ is based on an underdetermined problem. In that case we get $f_{2i}^{k+1} = f_i^k$ (see Appendix A.1), hence $S_n^{2n-1}$ is the interpolating Dubuc–Deslauriers $2n$-point scheme.

The following observation shows that it is sufficient to consider only odd degrees $d$.

**Proposition 7.** *For given data $y_1, \ldots, y_m$ at equidistant points $x_1, \ldots, x_m$ with $x_i = a + ih$, let $p$ and $q$ be the least squares polynomials of degrees $2k$ and $2k+1$, respectively. Then,*

$$p(c) = q(c),$$

*for $c = (x_1 + \cdots + x_m)/m$.*

The proof of this proposition is given in Appendix A.3. As a consequence of Proposition 7, the primal $2n$-point least squares schemes $S_n^{2k}$ and $S_n^{2k+1}$ are identical. This also means that the scheme of degree $2k$ reproduces polynomials of one degree more than expected by construction, matching the observation in [11] that the reproduction of odd degree polynomials comes "for free" by the primal symmetry.

We can also generalize the construction in Section 4.1 and define the dual $2n$-point least squares scheme of degree $d$ by the refinement rules

$$f_{2i}^{k+1} = \tilde{p}_{n,i}^d\big((3t_i^k + t_{i+1}^k)/4\big) \qquad \text{and} \qquad f_{2i+1}^{k+1} = \tilde{p}_{n,i}^d\big((t_i^k + 3t_{i+1}^k)/4\big),$$

which simplify to the rules in (16) for $d = 1$. Like $S_n^d$, this scheme $\bar{S}_n^d$ reproduces polynomials of degree $d$ by construction and its approximation order is $O(h^{d+1})$. Moreover, the scheme $\bar{S}_n^{2n-1}$ is the dual $2n$-point scheme [10].

Similar constructions lead to primal and dual $(2n+1)$-point least squares schemes of degree $d$, but we omit the details as they are straightforward. Apart from the increased approximation order, these schemes also tend to have a higher smoothness. For example, we verified numerically that the schemes $\bar{S}_n^3$ generate $C^3$ limit functions for $n = 4$ and $n = 5$, but we do not recommend to use them, because the rules become more complicated and the benefit of using them for reconstructing functions from noisy data is marginal, as shown in the next section.

## 4.4 Simulations of the primal least squares schemes of higher degree

The statistical model which is presented in Subsection 3.1 is valid for schemes based on higher degree least squares polynomials as well due to linearity. Furthermore, the first part of Theorem 6 is also true. However, the bounds of $\psi$ becomes puzzling due to the fact that the mask of $S_n^d$ with $d > 1$ is no longer positive nor given explicitly. Nevertheless, most of the conclusions can be seen through numerical trials. For example, in Table 3 we present several maximum and minimum values of $\psi_n$, as well as $\int_0^1 \psi_n(x)dx$ for $m > 1$. Thus, we conjecture

**Conjecture 8.** *Let $\psi_n$ be defined by (14) for schemes based on least squares polynomials of high degree. Then*

1. $\lim_{n \to \infty} (\max(\psi_n)) = 0$.

2. $\int_0^1 \psi_{n_1}(x)dx \geq \int_0^1 \psi_{n_2}(x)dx$, *for* $n_1 < n_2$.

This conjecture is also supported by the direct calculation of the basic limit functions, as in Figure 2, and by the simulations illustrated in Figure 3 and Figure 9. These figures present the function $\psi_n$ with different values of $m$ (the degree) and $n$ (the support size). In addition, they also evident the observation that as $n$ grows $\psi$ becomes almost constant.

The deterministic error in (13) is strongly related to $m$. This can be seen by the polynomial reproduction property of our schemes, that is the reconstruction of any $m$-th degree polynomial from its values at the integer by the limit of $S_n$. The latter property implies that the approximation order is (at least) $\mathcal{O}(h^{m+1})$. Thus, for larger $m$ the contribution of the deterministic error decreases, while we conjecture that the effect of the noise increases. In particular,

**Conjecture 9.** *Let $\psi_n$ be the function defined in (14) for $S_n$, a scheme based on least squares polynomials of degree $m_1$, and let $\hat{\psi}_n$ be the corresponding function for $\hat{S}_n$, a scheme based on least squares polynomials of degree $m_2$, where $m_1 \geq m_2$. Then,*

$$\int_0^1 \psi_n(x)dx \geq \int_0^1 \hat{\psi}_n(x)dx.$$

As show in Subsection 3.1, there is a trade-off between the deterministic approximation error and the effect of the noise on the expected square error. In particular, higher values of $n$ decrease the effect

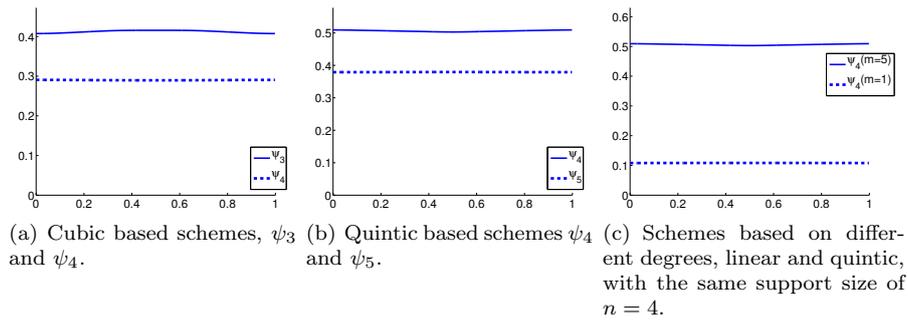

(a) Cubic based schemes, $\psi_3$ and $\psi_4$.  (b) Quintic based schemes $\psi_4$ and $\psi_5$.  (c) Schemes based on different degrees, linear and quintic, with the same support size of $n = 4$.

Figure 9: $\psi_n$ for different values of $m$ (the degree) and $n$ (the support size). Note the difference in scale on each sub-figure. We see that $\psi$ decreases where $n$ grows for $m > 1$ as well, and that $\psi$ increases where $n$ is fixed and $m$ grows (Figure 9(c)).

| $d, n$ | $\min(\psi_n)$ | $\max(\psi_n)$ | $\int_0^1 \psi_n(x)dx$ |
|---|---|---|---|
| 1,1 | .5 | 1 | 0.6647 |
| 1,3 | 0.1484 | 0.1489 | 0.1485 |
| 1,5 | 0.0847 | 0.0849 | 0.0847 |
| 1,7 | 0.0591 | 0.0592 | 0.0591 |
| 3,2 | 0.6406 | 1 | 0.7990 |
| 3,3 | 0.4074 | 0.4156 | 0.4115 |
| 3,5 | 0.2252 | 0.2254 | 0.2252 |
| 3,7 | 0.1563 | 0.1565 | 0.1564 |
| 5,3 | 0.7060 | 1 | 0.8447 |
| 5,5 | 0.3790 | 0.3793 | 0.3791 |
| 5,7 | 0.2573 | 0.2574 | 0.2573 |

Table 3: Several maximum, minimum, and average values of $\psi_n$, for different values of $m$ (the degree) and $n$ (the support size). The integrals are calculated using the trapezoidal numerical integration for equispaced points with $h = 0.002$ over $[0, 1]$.

of the noise but increase the deterministic error due to averaging of the values $\{f_i\}$ by weights with a large support $\phi_n(x - i)$, $|i| \leq 2n$. We believe that this behaviour occurs also for $m > 1$ on the ground of the second part of Conjecture 8. On the other hand, by increasing $m$ and keeping $n$ fixed increases the effect of the noise in view of Conjecture 9, while the deterministic error decreases.

We illustrate the schemes based on higher degree polynomials by two examples. We use two cubic based schemes ($m = 3$) with different supports, $n = 6$ and $n = 9$. Again, we present the limits of the subdivision schemes operate on data which is contaminated with noise. The sampled function is $f(x) = \cos(0.1x) - (\frac{x}{50} - 1)^3$. Figure 10 depicts this experiment.

In the last case, the smaller supported scheme with $n = 6$ responds more to the presence of the noise and thus approximates worse the true function. However, as seen in Figure 11, when we sampled function is the oscillatory function $f(x) = \cos(0.4x) - (\frac{x}{50} - .8)^3$, the more local scheme ($n = 6$) approximates $f$ better than the scheme $S_9$, since that for such a function the error is affected more by the deterministic error than by the noise.

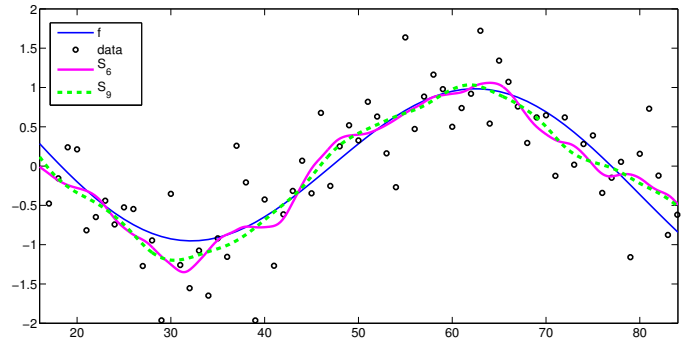

Figure 10: Two schemes based on cubic polynomials ($m = 3$) with different supports. High level of noise.

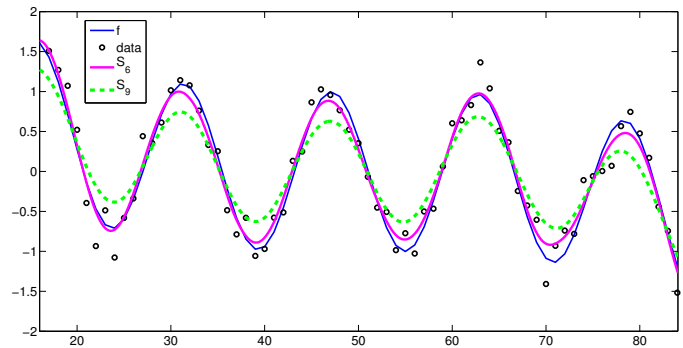

Figure 11: Two schemes based on cubic polynomials ($m = 3$) with different supports. An oscillatory sampled function.

**Acknowledgement.** We wish to thank Prof. Felix Abramovich for his help in formulating the statistical model for the expected square error.

# A  Least squares polynomials and least squares schemes

In this appendix we derive explicitly several properties of least squares polynomials, which were used throughout this paper. Some of the properties can be considered as common knowledge. However, for the completeness of the paper and in order to keep it as self-contained as possible, we present them here.

## A.1 Preliminaries and notations

Our subdivision schemes are based on least squares polynomial fitting. We denote by $\Pi_d$ the space of all polynomials of degree at most $d$. The fitting of data of the form $(x_i, y_i)$, $i = 0, \ldots, n$ by a polynomial of degree $d$, where $n \geq d$, depends on the $d + 1$ coefficients of the polynomial $p^* = \sum_{j=0}^{d} \beta_j x^j \in \Pi_d$, which satisfy

$$p^* = \arg\min_{p \in \Pi_d} \sum_{i=0}^{n} (p(x_i) - y_i)^2. \tag{19}$$

The coefficient $\boldsymbol{\beta} = (\beta_0, \ldots, \beta_d)$ of the polynomial $p^*$ are typically computed by differentiating the minimized functional in (19). This results in the normal equations

$$A^T A \boldsymbol{\beta} = A^T \boldsymbol{y},$$

where $A$ is an $(n+1) \times (d+1)$ matrix with entries $A_{i,j} = x_i^j$, $A^T$ is the transpose of $A$, and $\boldsymbol{y} = (y_0, \ldots, y_n)$ is the data vector. The matrix $A^T A$ is invertible for any set of distinct data points $x_i$, $i = 0, \ldots, n$, and the solution of the normal equations is given by

$$\boldsymbol{\beta} = \left(A^T A\right)^{-1} A^T \boldsymbol{y} = A^\dagger \boldsymbol{y}, \tag{20}$$

with $A^\dagger = \left(A^T A\right)^{-1} A^T$ the generalized inverse of $A$ or the "pseudo" inverse matrix. The matrix $A^\dagger$ is also known as the Moore-Penrose inverse [2], and is used to solve (19) directly from $A\boldsymbol{\beta} = \boldsymbol{y}$.

**Remark 10.** In case $d = n$, $p^*$ is the unique *interpolating* polynomial to the data. Furthermore, this ansatz can also be used in the case $d > n$ to pick among all interpolating polynomials the one with the smallest $\ell_2$-norm of coefficients $\|\boldsymbol{\beta}\|_2$, but then the solution depends on the particular basis of $\Pi_d$ that one chooses. However, $p^*(x_i) = y_i$, $i = 0, \ldots, n$, independently of the choice of the basis of $\Pi_d$.

Next we express the solution of (19) for $n \geq d$ in term of orthogonal polynomials. For that matter, we recall the notion of orthonormal polynomials relative to a discrete inner product. Let $\mathbf{X} = \{(x_i)\}_{i=0}^{n}$ be discrete equispaced data points, namely $x_i = x_0 + ih$, $i = 0, \ldots, n$. We denote by $\langle, \rangle$ the standard inner product in Euclidean space. In addition, for functions $f, g : \mathbf{X} \to \mathbb{R}$ we define the reduced inner product

$$\langle f, g \rangle_{\mathbf{X}} = \langle f|_{\mathbf{X}}, g|_{\mathbf{X}} \rangle, \tag{21}$$

where $f|_{\mathbf{X}}$ is the restriction of $f$ to the set $\mathbf{X}$. A family $\mathcal{L} = \{L_0(x), L_1(x), \ldots, L_k(x)\}$ of $k$ polynomials in $\Pi_d$ with $d \geq k$, is orthonormal with respect to the inner product (21) if

$$\langle L_i, L_j \rangle_{\mathbf{X}} = \delta_{i,j}, \quad i, j = 0, \ldots, k, \tag{22}$$

where $\delta_{i,j}$ is the standard Kronecker delta, namely $\delta_{i,j} = 1$ if $i = j$ and zero otherwise. Furthermore, under the assumptions $L_j \in \Pi_j$, $j = 0, \ldots, k$ there exists a unique family $\mathcal{L}$ satisfying (22). For more details see [13] and reference therein.

The use of such orthonormal polynomials for the least squares problem is well-known. In particular, $p^* \in \Pi_m$ of (19) has the form [13]

$$p^*(x) = \sum_{i=0}^{m} \langle L_i|_{\mathbf{X}}, \boldsymbol{y} \rangle L_i(x). \tag{23}$$

## A.2 Derivation of the masks

A naive implementation of the refinement rules (16) is computationally expensive, since the solution of the least squares problem is equivalent to the solution of a linear system. Using (23) we can exploit the shift invariance and scale invariance properties of the least squares polynomials, which guarantee a much simpler way for the evaluation of the refinement rules. This is shown in the following proposition.

**Proposition 11.** *Let $\boldsymbol{f}^k = \{f_i^k\}_{i \in \mathbb{Z}}$ be the data of the $k$-th level of our subdivision. The refinement rules (18) of $S_n^d$ are independent of $k$, and have the form*

$$f_i^{k+1} = \left(S_n^d(\boldsymbol{f}^k)\right)_i = \sum_j \alpha_{i-2j} f_j^k, \tag{24}$$

*where only a finite number of the coefficients $\{\alpha_i\}_{i \in \mathbb{Z}}$ are non-zero in (24).*

*Proof.* Consider (23). By using the definition of the inner product (21) we have

$$p^*(x) = \sum_{i=0}^n \ell_i(x) y_i, \tag{25}$$

with $\ell_i(x)$ the polynomial

$$\ell_i(x) = \sum_{j=0}^d L_j(x_i) L_j(x) \tag{26}$$

Since in each least squares subdivision scheme, the degree of the least squares polynomials determining the refinement rule is fixed, it is sufficient to show that $\ell_i(x)$ does not depend on $h$ or $x_0$. Indeed, the orthonormal polynomials are invariant under affine transformations of the points $\mathbf{X}$. To be specific, let $\mathbf{X} = \{x_0, x_0 + h, \ldots, x_0 + nh\}$, $n > d$, and let $\{L(x)\}_{i=0}^d$ denote the corresponding orthonormal polynomials with respect to $\langle \cdot, \cdot \rangle_{\mathbf{X}}$. For

$$\bar{\mathbf{X}} = \{\bar{x}_0, \bar{x}_0 + \bar{h}, \ldots, \bar{x}_0 + n\bar{h}\},$$

define

$$\bar{L}_i(x) = L_i(x_0 + \frac{x - \bar{x}_0}{\bar{h}} h).$$

Then

$$\bar{L}_i(x)(\bar{x}_0 + \lambda \bar{h}) = L_i(x_0 + \lambda h), \quad \lambda \in \mathbb{R}. \tag{27}$$

The latter implies that

$$\langle \bar{L}_i(x), \bar{L}_j(x) \rangle_{\bar{\mathbf{X}}} = \sum_{s=0}^n \bar{L}_i(x)(\bar{x}_0 + s\bar{h}) \bar{L}_j(x)(\bar{x}_0 + s\bar{h}) = \langle L_i, L_j \rangle_{\mathbf{X}} = \delta_{i,j}.$$

Thus, the set $\{\bar{L}_i(x)\}$ is a family of orthonormal polynomials over $\bar{\mathbf{X}}$. In view of (26) and (27) $\ell_i(x)$ is independent of $h$ and $x_0$. □

The finite set of coefficients $\{\alpha_j\}$ of Proposition 11 is called the *mask* of $S_n^d$. As a direct consequence of Proposition 11,

**Corollary 12.** *The least squares subdivision schemes are stationary linear schemes.*

## A.3 Proof of Proposition 7

The form (23) of the leasts squares orthogonal polynomial leads to the following proof of Proposition 7.

*Proof.* The refinement rules (18) evaluate the least squares polynomials at the center of their data points. Due to the fact that the points are equidistant we know that $L_{2i}$ is an even polynomial about the center of its data points and $L_{2i+1}$ is an odd polynomial about that center [13]. Thus, the value of each $L_{2i+1}$ is zero at the center of the data points. Therefore, by (25) and (26) the $2d$ and $2d+1$ least squares polynomials coincide at this center, which concludes the proof. □

## A.4 Computation of the masks

The explicit form of (2), for the schemes based on linear least squares polynomials can be easily concluded from (23). The linear least squares polynomial is of the form

$$p^*(x) = \langle L_0|_{\mathbf{X}}, \mathbf{y}\rangle L_0(x) + \langle L_1|_{\mathbf{X}}, \mathbf{y}\rangle L_1(x).$$

The polynomial $p^*$ is evaluated at the centre of the data points $\mathbf{X}$. As noted in the proof of Proposition 7, $L_1(x)$ vanishes at this center point. It is straightforward to see that $L_0 = \frac{1}{\sqrt{n+1}}$. Hence, for the center point $c$ we have

$$p^*(c) = \langle L_0|_{\mathbf{X}}, \mathbf{y}\rangle L_0(c) = \frac{1}{n+1}\sum_{j=0}^n y_j.$$

We present now an efficient algorithm for the implementation of our least square schemes of high degree. Proposition 11 suggests a method for calculating the masks based on (25) and (26) (explicit form of the orthogonal polynomials can be found in [13]). Here we introduce an alternative method for calculating $\ell_i$, $i = 0, \ldots, n$ in (25), which follows from the linearity of $p^*$ in $\mathbf{y}$, see (20).

**Proposition 13.** *In* (25) $\ell_i(x)$ *is the least squares polynomial*

$$\ell_i(x) = \arg\min_{p\in\Pi_d}\sum_{j=0}^n \left(p(j) - \delta_{i,j}\right)^2,$$

*where $\delta_{i,j}$ is the standard Kronecker delta.*

We use the form (25) of the least squares polynomial with $\ell_i$ as in Proposition 13 to calculate the mask of $S_n^d$ as follows. We start with the calculation of two "pseudo" inverse matrices, as in (20), over the data points $-n+1, \ldots, n-1$ and $-n+1, \ldots, n$. The first data points correspond to the refinement rule of the even indices of (18) while the second data points correspond to the refinement rule of the odd indices of (18). Note that by the proof of Proposition 11 we can also use a different scale of points, that is $(-n+1)h, \ldots, (n-1)h$ and $(-n+1)h, \ldots, nh$ with a small $h < 1$, in order to have more stable pseudo inverse matrices.

To calculate the mask we next determine $\ell_i$ of Proposition 13 for each of the refinement rules (18). Note that in this point, for the coefficients of $\ell_i$ we only have to multiply a $4n-1$ length vector with a pseudo inverse matrix. Finally we obtain the mask by evaluating $p^*$ by the polynomials $\ell_i$ of the first even refinement at 0, and those of the second odd refinement at $\frac{1}{2}$. The calculations of the masks of the dual schemes are done similarly.